\documentclass[11pt,amssymb]{article}
\usepackage{amsfonts,amsmath}

\def \beq {\begin{eqnarray}}
\def \eeq {\end{eqnarray}}
\def \beqn {\begin{eqnarray*}}
\def \eeqn {\end{eqnarray*}}

\newcommand{\halmos}{\rule{1ex}{1.4ex}}

\newcounter{for}[section]

\numberwithin{equation}{section}

\newtheorem{itlemma}{Lemma}[section]
\newtheorem{itproposition}[itlemma]{Proposition}
\newtheorem{theorem}[itlemma]{Theorem}
\newtheorem{itcorollary}[itlemma]{Corollary}
\newtheorem{itremark}[itlemma]{Remark}
\newtheorem{itremarks}[itlemma]{Remarks}
\newtheorem{itdefinition}[itlemma]{Definition}
\newtheorem{itexample}[itlemma]{Example}

\newenvironment{fact}{\begin{itfact}\rm}{\end{itfact}}
\newenvironment{claim}{\begin{itclaim}\rm}{\end{itclaim}}
\newenvironment{lemma}{\begin{itlemma}}{\end{itlemma}}
\newenvironment{remark}{\begin{itremark}\rm}{\end{itremark}}
\newenvironment{remarks}{\begin{itremarks} \rm}{\end{itremarks}}
\newenvironment{corollary}{\begin{itcorollary}}{\end{itcorollary}}
\newenvironment{proposition}{\begin{itproposition}}{\end{itproposition}}
\newenvironment{definition}{\begin{itdefinition}\rm}{\end{itdefinition}}
\newenvironment{example}{\begin{itexample}\rm}{\end{itexample}}
\newenvironment{proof}{\noindent {\em Proof}.\ \
}{\hspace*{\fill}$\halmos$\medskip}
\newcommand{\be}[1]{\addtocounter{for}{1} \begin{equation}\label{#1}}
\newcommand{\ee}{\end{equation}}
\newcommand{\bl}[1]{\begin{lemma}\label{#1}}
\newcommand{\br}[1]{\begin{remark}\label{#1}}
\newcommand{\brs}[1]{\begin{remarks}\label{#1}}
\newcommand{\bt}[1]{\begin{theorem}\label{#1}}
\newcommand{\bd}[1]{\begin{definition}\label{#1}}
\newcommand{\bp}[1]{\begin{proposition}\label{#1}}
\newcommand{\bc}[1]{\begin{corollary}\label{#1}}
\newcommand{\bfact}[1]{\begin{fact}\label{#1}}
\newcommand{\bex}[1]{\begin{example}\label{#1}}
\newcommand{\ec}{\end{corollary}}
\newcommand{\efact}{\end{fact}}
\newcommand{\eex}{\end{example}}
\newcommand{\el}{\end{lemma}}
\newcommand{\er}{\end{remark}}
\newcommand{\ers}{\end{remarks}}
\newcommand{\et}{\end{theorem}}
\newcommand{\ed}{\end{definition}}
\newcommand{\ep}{\end{proposition}}
\newcommand{\epr}{\end{proof}}
\newcommand{\bpr}{\begin{proof}}
\newcommand{\bcl}[1]{\begin{claim}\label{#1}}
\newcommand{\ecl}{\end{claim}}

\newcommand{\ecs}{\end{corollary}}
\newcommand{\eers}{\end{exercise}}
\newcommand{\eexs}{\end{example}}
\newcommand{\eems}{\end{example}}
\newcommand{\els}{\end{lemma}}
\newcommand{\eles}{\end{lemmaex}}
\newcommand{\ets}{\end{theorem}}
\newcommand{\eds}{\end{definition}}
\newcommand{\eps}{\end{proposition}}

\def\D{{\mathbb D}}
\def\N{{\mathbb N}}
\def\Z{{\mathbb Z}}

\def\NN{\mathcal N}
\def \ind {\hbox{1\hskip -3pt I}}

\def\bs{\backslash}
\def\reff#1{(\ref{#1})}
\newcommand {\cro}[1] {\left[ {#1} \right]}
\newcommand {\acc}[1] {\left\{ {#1} \right\}}
\newcommand {\pare}[1] {\left( {#1} \right)}

\begin{document}

\title{A note on the rightmost particle\\
in a Fleming-Viot process.} 
\author{Amine Asselah \&  Marie-No\'emie Thai\\
\sl Universit\'e Paris-Est}
\maketitle
\noindent {\bf Abstract}
We consider $N$ nearest neighbor random walks on the positive
integers with a drift towards the origin. When one 
walk reaches the origin, it jumps to the position of
one of the other $N-1$ 
walks, chosen uniformly at random. We show that this particle
system is ergodic, and establish some exponential moments
of the rightmost position, under the stationary measure.

\vskip 3mm



\section{Introduction}
There is recent interest in approximating the limiting law
of irreducible Markov processes conditioned not to hit some
(forbidden) state \cite{Lelievre,Villemonais}. This limiting law
is not guaranteed to exist, but when it does it is called a 
quasi-stationary distribution (QSD). For QSD in
the context of Birth and Death chains, we refer to \cite{Ferrari-Kesten},
and the situation treated there is one in which there is a one
parameter family of QSD.

QSD are neither well understood, nor easily amenable
to simulation. One proposal made by Burdzy, Holyst, Ingerman, and
March~\cite{Burdzy} (in a particular setting)
is to consider $N$ independent Markov processes except that when
one reaches the forbidden state, it jumps to the state of one of the other
processes, chosen uniformly at random. The natural conjecture
is that the empirical measure, under the stationary
measure, converges to the QSD as the number $N$ of processes
goes to infinity. It is also natural to
conjecture that the selected QSD is {\it the minimal}, in terms
of average time needed to reach the forbidden state.

In this note, we consider $N$ random walks on $\N$, 
in continuous time, with a drift
towards the origin. When one random walk reaches
the origin, it jumps instantly to the position of one of
the other $N-1$ walks, choosen uniformly at random.
We call Fleming-Viot the interacting random walks just described.
Indeed, this dynamics has a genetic interpretation.
\begin{itemize}
\item
Positions of the walks are evolving genetic traits of $N$ individuals.
\item
The forbidden state (here 0) is a lethal trait ({\it the selection
mechanism}). 
\item At the moment an individual dies, another
one, chosen uniformly at random, branches ({\it the branching
mechanism}). This keeps the population size constant.
\end{itemize}

We establish a Foster's criteria, which gives
ergodicity, as well as a control of
small exponential moments of the rightmost walk.
To state our main result,
let $\xi_T$ denote the position of the $N$ interacting walks 
at time $T$, and let $E[\cdot|\xi_0=\xi]$ denote 
average with respect of the law of the process $\xi_T$
with initial condition $\xi$.

\bt{theo-1} 
There are positive constants $K,\alpha,\kappa,
\delta_0,A,c_1,c_2,c_3$ such that
for $N\in \N$, time $T=A\log(N)$, any $\delta<\delta_0$,
and $\xi\in\N^N$, we have
\be{foster}
\begin{split}
E\cro{\exp\big(\delta\max(\xi_T)\big)\big|\xi_0=\xi}-&\exp\big(\delta
\max(\xi)\big)<
-\,c_1\, \ind_{\max(\xi)>K\log(N)}e^{\delta\max(\xi)}\\
& + c_2 \ind_{\max(\xi)>K\log(N)} e^{-\kappa T} e^{\delta\max(\xi)}+
c_3e^{\delta\alpha \log(N)}.
\end{split}
\ee
As a consequence, for $N$ large enough there is a unique invariant
measure $\lambda_N$ for Fleming-Viot.
Integrating over $\lambda_N$, there are $\beta,C>0$ such that for any $N$,
and $\delta<\delta_0$
\be{rightmost-expon}
\int \exp(\delta\max (\xi)) d\lambda_N(\xi) \;\le\; 
C\exp\big(\delta \beta\log(N)\big).
\ee
\et

This first elementary step is an important ingredient 
in the proof of the conjecture we alluded to above \cite{Thai}.
Also, it might be of independent interest in view of recent deep and
comprehensive studies on the rightmost position in branching random walks
\cite{Derrida-Brunet,Addario-Reed,
Hu-Shi,Bovier1,Aidekon-Shi,Bovier2,Aidekon,Maillard}.
This selection of recent works is far from being exhaustive, 
but already shows the vitality of this issue.

The rest of paper is organized as follows. In Section~\ref{sec-model}
we define the model, and recall well-known large deviations estimates.
In Section~\ref{sec-independence}, we explain how
to divide walks into groups with
little correlations over a well chosen time period. 
In Section~\ref{sec-bad}, we estimate the probability that the
maximum displacement does not decrease.
Finally, in Section~\ref{sec-foster}, we establish Foster's criteria.

\section{Model and Preliminaries}\label{sec-model}
Here, we deal with continuous-time nearest neighbor random walks on $\N$, 
with rate $p$ to jump right, and rate $q=1-p>p$ to jump left. The
drift is $-v$ with $v=q-p>0$. A single walk makes $N_t$ jumps
in the time period $[0,t]$, and its increments are denoted
$X_1,\dots,X_n$, with $E[X_i]=-v$, and $\bar X_i=X_i+v$ denotes
the centered variable. Note that
\be{ld-1}
P\big(\sum_{i=1}^T(X_i+v)\ge  xT\big)\le \exp(-T I(x)),
\ee
with
\be{def-I}
I(x)=\sup_{\lambda>0}\acc{ \lambda x-\Lambda(\lambda)}\quad
\text{with}\quad \Lambda(\lambda)=\log(pe^\lambda+qe^{-\lambda})+\lambda v.
\ee
Due to the nearest neighbor jumps of our walk, we have
\be{cond-I}
I(v+1)=\log\big(\frac{1}{1-q}\big)\quad\text{and for}\quad
x>v+1,\quad I(x)=\infty.
\ee
We define also $x\mapsto \tilde I(x)=1-\exp(-I(x))$, which is
discontinuous at $v+1$ with
\be{cond-tildeI}
\tilde I(v)<\tilde I(v+1)=q\quad\text{and for}\quad x>v+1,\quad \tilde I(x)=1.
\ee
Note that if $N_T$ is Poisson of mean $T$, then
\be{ld-2}
P\big(\sum_{i\le N_T}(X_i+v)\ge  xN_T\big)\le \exp(-T \tilde I(x)).
\ee

\paragraph{On Poisson tails.} We need two rough tail estimates
on the Poisson clocks. Both are obvious and well-known. Assume $\chi,T$
are positive.
\be{poisson-ut}
P(N_T\ge eT+\chi)\le \exp(-T-\chi),
\ee
and
\be{poisson-lt}
P(N_T\le \frac{1}{e}T-\chi)\le \exp(-(1-2/e)T-\chi),
\ee
Both are obtained readily by Chebychev's inequality. Indeed, we
obtain \reff{poisson-ut} from
\be{step-ut}
P(N_T\ge eT+\chi)\le e^{-eT-\chi}E[e^{N_T}]=\exp(-T-\chi),
\ee
and we obtain \reff{poisson-lt} from
\be{step-lt}
P(N_T\le \frac{1}{e}T-\chi)\le e^{T/e-\chi}E[e^{-N_T}]=\exp(-(1-2/e)T-\chi).
\ee

\section{Independence}\label{sec-independence}
\paragraph{On the multitype branching of \cite{AFGJ}.}
A key idea introduced in \cite{AFGJ} is to embed the Fleming-Viot
process into a multitype branching process whose space displacements
and branching mechanism are independent, and which is {\it attractive}. 
We refer to Section 3 of \cite{AFGJ} for a description of the multitype
branching process, and recall here its main features.
Assume that we start with $N$ interacting random walks.
This defines $N$ types with which we associate $N$ 
independent exponential clocks of intensity $q$ with marks.
The time realizations of the
clock of type $i$ have marks in the set of labels $\{1,\dots,
N\}\bs\{i\}$, and each mark is chosen uniformly at random from
the $N-1$ symbols. When clock $i$ rings, and when its mark is $j$,
each walk of type $j$ branches into two children: one of type
$i$ and one of type $j$. The two children behave as independent random
walks starting at the position of their parent. If $\D_T$ denotes
the population of individuals alive at time $T$, and $|\D_T|$ denotes
its cardinal, it is easy to see the equality $E[|\D_T|]=|\D_0|\exp(q T)$.
For an individual $v$ alive at time $T$, we denote by $t\mapsto S_v(t)$
its trajectory for $t\in [0,T]$.

\paragraph{Independent groups of walks.}
A drawback of the multitype branching process
is an exponentially growing population. Since, we use
a time of order $\log(N)$, we cannot use here such an embedding.
Even though in the Fleming-Viot process, all particles
interact which each other, a simple observation is that 
as long as a particle has not touched the origin its
trajectory is independent from the other ones,
even though this trajectory might influence others.

To create some independence between walks, 
we decompose the interacting walks in two sets at time 0. We first
fix a time $T$ and a length $L$ to be chosen later.
\begin{itemize}
\item The {\it blacks}, whose initial position is below $L$.
\item The {\it reds}, whose initial position is above $L$.
\end{itemize}
Then, color changes as follows: if a black walk jumps
on a red walk, it becomes instantly red. We interpret this
jump as a {\it red binary branching}. Now, red walks are
not independent from black walks because they might touch
the origin before time $T$, and jump onto a black position.
However, if $vT\ll L$, we expect this to be rare. To obtain
independence, we add another
color to our description: each red walk is coupled
with a green walk which behaves identically in terms of move
or branching but with green children,
except that when a green walk reaches the origin
it continues its drifted motion on $\Z$ (without selection mechanism).
Thus, green walks behaves like independent random walks with
branching at the times a black particle hit zero and chooses
the label of a green walk. If $R_0$ is the first time
one of the red walks touches the origin, we have that at time $T$,
on the event $\{R_0>T\}$,
red and green positions are identical. The point of introducing
green walks is that their branching times is independent of
their positions. We denote with $D_T^r,D_T^g,D_T^b$ the 
respective number of red, green and black walks at time $T$.
Also $D_T=D_T^r\cup D_T^g\cup D_T^b=\{1,\dots,N\}$, and we still denote
by $t\mapsto S_v(t)$, the trajectory of $v\in D_T$.

When embedding a group of walks into a branching process, we denote with
$\D_T^r,\D_T^g,\D_T^b$ the respective number of red, green and black 
individuals in the mutlitype branching processes. 

The key idea here is to work on a time of order $\log(N)$,
to control the black walks by a multitype branching process,
but to let the red walks (or rather the green walks)
grow as in Fleming-Viot with a population bounded by $N$, and with
branching due to independent black walks.

\paragraph{On the choice of time $T$ and length $L$.}
We choose $T$ large enough so that $q+\log(N)/T<1$.
We actually need a little more.
\be{def-T}
\kappa=\min\big(1-2/e,\tilde I(\frac{v}{2})-\frac{\log(N)}{T},
1-q-\frac{\log(N)}{T}\big)>0.
\ee
Once $T=A\log(N)$ satisfies \reff{def-T}, we set $L=eT$.

\section{When things go wrong}\label{sec-bad}
We wish to estimate the probability of the event where
the maximum displacement does not decrease.
We thus define a {\it bad set} $B(T,L)$ as containing the following
events:
\begin{itemize}
\item One red walk reaches the origin before time $T$ (i.e. $\{R_0\le T\}$).
\item One black walk travels a distance $L$ upwards in a period $[0,T]$.
\item The maximum displacement of a green walk in a time $T$ is
above $-\frac{v}{2e} T$.
\end{itemize}
Thus, on the complement on $B(T,L)$, green and red are identical,
and 
\[
\max_{v\in D_T} S_v(T)-\max_{v\in D_0} S_v(0)\le
\max_{v\in D_T}\Big( S_v(T)-S_v(0)\Big)
=\max_{v\in D_T^g}\Big( S_v(T)-S_v(0)\Big)< -\frac{v}{2e} T,
\]
which implies that if $M(T)=\max_{v\in D_T} S_v(T)$ 
\be{good-main}
E\cro{\ind_{B^c(T,L)}\exp\big(\delta(M(T)-M(0))\big)\big|\xi(0)=\xi}\le 
\exp\big(-\frac{v\delta}{2e} T\big).
\ee
We estimate next the probability of each event making up $B(T,L)$,
with the following outcome.
\bl{lem-bad}
For any $\xi\in \N^N$, we have, with $\kappa$ as in \reff{def-T},
\be{bad-1}
P(B(T,L)\big|\xi(0)=\xi)\le 4\exp(-\kappa T).
\ee
\el

\subsection{A red walk does reach 0}\label{sec-red}
Recall that $L=e T$. 
We embed the Fleming-Viot into a branching multitype, while
keeping the red coloring. We need
to estimate the probability that one red displacement 
gets below $L$ units in a time period $[0,T]$. Note that
to realize $\{R_0<T\}$, there is $v\in \D_T^r$ such that
the number of its time jumps $N_T$ must
be larger than $L$, and this is what we use.
\be{ori-1}
\begin{split}
P(R_0<T\big|\xi(0)=\xi)\le&
E[|\D^r_T|]\times P\big( \exists t\le T,\ \sum_{i\le N_t} X_i<
-eT\big)\\
\le &E[|\D^r_T|]\times P\big( N_T> eT)
\le Ne^{qT} e^{-T}\le e^{-\kappa T}.
\end{split}
\ee
\subsection{A large black displacement}\label{sec-black}
Recall that at the time a black reaches 0, and jumps on
a red walk, it ceases to be black to become red. 
We bound here the black walks with
a multitype branching, assuming that blacks do
only jump on blacks, with the effect that we
are overestimating the black population. The
estimates are similar to these of Section~\ref{sec-red}.
We use that to make $L$ steps right, a black walk must make
$L$ time-marks ($N_T>L$), and this event is estimated in \reff{ori-1}.

\subsection{Green's maximum too high}\label{sec-green}
The key point is that the green branching times are
independent of positions of the green. They depend
only on the history of black walks. Also, the
population of green walks is bounded by $N$. Thus, it
is crucial here not to use the multitype branching of \cite{AFGJ}:
we estimate the probability
that $\{\max_{v\in D_T^g}(S_v(T)-S_v(0))>- vT/(2e)\}$.
Define $\NN_T(\gamma)$ as the number of green walks whose
displacement during time period $[0,T]$ is larger than $\gamma$.
Then, 
\be{green-key}
E[\NN_T(\gamma)\big|\xi(0)=\xi]=E[|D_T^g| \big| \xi(0)=\xi]\times 
P\big( \sum_{i\le N_T} X_i> \gamma \big).
\ee
The reason is the independence of the branching times 
and displacements of the green walks.
For $v\in D_T^g$, assume there are $\nu$ branchings before
time $T$, say at times $T_1,\dots,T_\nu$, and we have 
(for $X_i^{(k)}$ i.i.d.  independent from $\{T_i, i\in \N\})$
\be{green-walk}
S_v(T)-S_v(0)=\sum_{i\in N[0,T_1]}X_i^{(1)}+\dots
+\sum_{i\in N[T_\nu,T]}X_i^{(\nu)}.
\ee
As one conditions first on the black history up to time $T$,
one fixes the times $T_1,\dots,T_\nu$, and obtain that
$N[0,T_1]+\dots+N[T_\nu,T]$ sums up to a Poisson variable $N[0,T]$
of intensity $T$, and most importantly
\be{green-7}
\sum_{i\in N[0,T_1]}X_i^{(1)}+\dots
+\sum_{i\in N[T_\nu,T]}X_i^{(\nu)}=\sum_{i\in N[0,T]} X_i,
\ee
where the $\{X_i,\ i\in \N\}$ are i.i.d. increments
independent of $N[0,T]$. We obtain, with $\kappa$ defined in \reff{def-T},
\be{green-1}
\begin{split}
P\big(\max_{v\in D_T^g}(S_v(T)-S_v(0))>&-
\frac{v}{2e}T\big|\xi(0)=\xi\big)\le 
E\cro{\NN_T\big(-\frac{vT}{2e}\big)\big|\xi(0)=\xi}\\
\le & N P\big( \sum_{i\le N_T} \bar X_i> vN_T-\frac{v}{2e}T\big)\\
\le & N\Big( P\big( \sum_{i\le N_T} \bar X_i> \frac{v}{2}N_T\big)+
P\big(N_T<\frac{1}{e} T\big)\Big)\\
\le & N \Big(\exp\big(-T \tilde I(\frac{v}{2})\big)+\exp\big(
-(1-2/e)T\big)\Big)\le 2 \exp(-\kappa T).
\end{split}
\ee

\section{Foster's criteria}\label{sec-foster}
We start with an estimate on the tail, and of the exponential moments.
\subsection{On exponential moments}\label{sec-moment}
We deal here with the multitype branching process.
Recall that $\D_0=\{1,\dots,N\}$, and let $S(0)=\{S_v(0),v\in \D_0\}$.
\bl{lem-tail}
For any $T$ satisfying \reff{def-T}, and any $\chi>0$
\be{exp-2}
P\big( \max_{v\in \D_T}\big( S_v(T)-S_v(0)\big) >e T
 +\chi\big|S(0)=\xi\big)\le \exp(-\chi).
\ee
\el
\bpr
Since the branching mechanism is independent of positions
\be{exp-3}
P\big( \max_{v\in \D_T}\big( S_v(T)-S_v(0)\big)
>eT+\chi\big|S(0)=\xi\big)\le
E[|\D_T|]\times P\big( \sum_{i=1}^{N_T} X_i>e T +\chi\big).
\ee
Now, if $\bar X_i$ denotes the centered variable, note that
since the walk is nearest neighbor
\[
P\big( \sum_{i=1}^{N_T}\bar X_i> (v+1)N_T)=0.
\]
Now,
\be{exp-4}
\begin{split}
P\big( \sum_{i=1}^{N_T} X_i>e T +\chi\big)=&
P\big( \sum_{i=1}^{N_T}\bar X_i>v N_T+e T +\chi\big)\\
\le & P\big( \sum_{i=1}^{N_T} \bar X_i>(v+1)N_T \big)+
P(N_T> eT+\chi)\\
\le & 0+P(N_T> eT+\chi)\\
\end{split}
\ee
Now, $N_T$ is a Poisson variable of mean $T$,
the standard estimate \reff{poisson-ut} leads to
\be{exp-5}
P\big(N_T> e T+\chi\big)\le \exp\big(-T-\chi\big).
\ee
Also, we have $E[|\D_T|]\le N\exp(qT)$, and \reff{exp-3} and
the choice of $T$ in \reff{def-T} yield
\be{exp-6}
P\big( \max_{v\in \D_T}\big( S_v(T)-S_v(0)\big) >e T +\chi
\big|S(0)=\xi\big)\le Ne^{qT}e^{-T-\chi}\le e^{-\chi}.
\ee
\epr

We can state our main estimate.
\bl{lem-exp}
Assume that $(1-q)T>log(N)$, and $\delta<1$. Then, we have
\be{exp-7}
E\cro{\exp\pare{\delta\big(\max_{i\le N} \xi_i(T)-\max_{i\le N} \xi_i(0)
\big)}\Big|\xi(0)=\xi}\le \frac{1}{1-\delta} e^{\delta e T}.
\ee
\el
\bpr
For any random variable $X$, we have
\be{exp-8}
\begin{split}
E\cro{e^{\delta X}}= &1+\int_0^\infty \delta e^{\delta u} 
P(X>u) du\\
\le & 1+\int_0^{e T} \delta e^{\delta u} du+
\int_{e T}^\infty \delta e^{\delta u} P(X>u) du\\
\le & e^{\delta e T}\Big( 1+\int_0^\infty \delta e^{\delta u}
P(X>u+e T) du\Big).
\end{split}
\ee
Now, using the tail estimate \reff{exp-2}, we have
\be{exp-9}
E\cro{\exp\pare{\delta\max_{v\in \D_T}
\big(S_v(T)-S_v(0)\big)}}\le e^{\delta e T}
\Big( 1+\int_0^\infty \delta e^{\delta u} e^{-u}\Big)
\le \frac{e^{\delta e T}}{1-\delta}.
\ee
We now use the following bound to conclude
\be{exp-0}
\begin{split}
E\cro{\exp\pare{\delta\big(\max_{i\le N} \xi_i(T)-\max_{i\le N} \xi_i(0)
\big)}\big|\xi(0)=\xi}\le&
E\cro{\exp\pare{\delta\big(\max_{v\in \D_T} S_v(T)-
\max_{v\in \D_0}S_v(0)\big)}\big|S(0)=\xi}\\
\le & E\cro{\exp\pare{\delta\max_{v\in \D_T}
\big(S_v(T)-S_v(0)\big)}\big|S(0)=\xi}.
\end{split}
\ee
\epr
\subsection{\bf Proof of Theorem~\ref{theo-1}}\label{sec-afgj}
We recall the general strategy
of the proof of Proposition 1.2 of \cite{AFGJ} (The Foster criteria).
We have a bad set $B(T,L)$ (which depends on $T$ and $L$)
which contains the cases where
the maximum increases over a period $[0,T]$, or when
black walks win over or influence red ones.
First, there is a set $K$ on which we do not expect the maximum
to decrease, with
\be{foster-3}
K=\acc{\max_v (S_v(0))< 3L}.
\ee
Then, there is a good set where the maximum decreases:
\be{foster-2}
K^c\cap B^c(T,L)\subset 
G=\acc{\max_{v\in D_T}(S_v(T))- \max_{v\in D_T}(S_v(0))\le -\frac{v}{2e} T}.
\ee
Now, set $M_t=\max S_v(t)$.
When $\xi$ is the initial configuration,
and when we work with $\xi\in K^c$, we have using Cauchy-Schwarz
\be{foster-4}
\begin{split}
\ind_{\xi\in K^c}&\Big(E\cro{e^{\delta M_T}\big|\xi(0)=\xi}-e^{\delta M_0}
\Big) =\ind_{\xi\in K^c}e^{\delta M_0}
\pare{E\cro{e^{\delta (M_T-M_0)}
\big(\ind_{B(T,L)}+\ind_G\big)\big|\xi(0)=\xi}-1}\\
\le &\ind_{\xi\in K^c}e^{\delta M_0} \Big(P\big(B(T,L)\big|\xi(0)=\xi\big) 
E\cro{\exp(2\delta(M_T-M_0))\big|\xi(0)=\xi}\Big)^{1/2}\\
&\qquad -\ind_{\xi\in K^c}e^{\delta M_0}\Big(1- e^{-\delta vT/2}\Big).
\end{split}
\ee
We know from Lemma~\ref{lem-exp} that for $\delta<1/2$ 
\[
E\cro{\exp(2\delta(M_T-M_0))\big|\xi(0)=\xi}\le 
\frac{e^{2\delta eT}}{1-2\delta}.
\]
Note also that on the set $K$, if we use Lemma~\ref{lem-exp}
\be{foster-5}
\ind_K\Big(E\cro{ e^{\delta M_T}\big|\xi(0)=\xi}- e^{\delta M_0}\Big)\le
\ind_K \exp( 3 \delta L+\delta eT).
\ee
Thus, adding \reff{foster-4} and \reff{foster-5}, we obtain
\be{foster-6}
E\cro{ e^{\delta M_T}}- e^{\delta M_0}\le
\ind_K e^{3\delta L+\delta eT}-\ind_{K^c}
\Big(1- e^{-\delta vT/2}\Big)e^{\delta M_0}+
\ind_{K^c} \Big(P(B(T,L)\big|\xi(0)=\xi)  
\frac{e^{2\delta eT}}{1-2\delta}\Big)^{1/2} e^{\delta M_0}.
\ee
We use now Lemma~\ref{lem-bad}, and choose $\delta$ small
enough so that $\kappa>4 \delta e$ with the result
\be{foster-8}
E\cro{ e^{\delta M_T}}- e^{\delta M_0}\le
e^{3\delta L+\delta eT} -\ind_{K^c}
\Big(1- e^{-\delta vT/2}\Big)e^{\delta M_0}+
\ind_{K^c}\frac{e^{-\kappa T/4}}{\sqrt{1-2\delta}}e^{\delta M_0}.
\ee
Inequality \reff{foster-8} is a
Foster's criteria (see \cite[Theorems 8.6 and 8.13]{Robert}).
This implies the first part of Theorem~\ref{theo-1}.

Now, as we integrate \reff{foster-8} with respect 
to the invariant measure, the
left hand side of \reff{foster-6} vanishes, and we obtain
\be{foster-7}
\Big(1- e^{-\delta vT/2}\Big) \int_{K^c} e^{\delta M(\xi)} d\lambda^N(\xi)
\le \exp( 3 \delta L+\delta eT)+ \exp(-\frac{\kappa}{4} T)
\int_{K^c} e^{\delta M(\xi)} d\lambda^N(\xi).
\ee
With $A$ large enough so that \reff{def-T} holds
with $T=A\log(N)$ and $L=eT$, the second part of
Theorem~\ref{theo-1} follows at once.

\paragraph{\bf Acknowledgements.}
We would like to thank Djalil Chafai for valuable discussions.

\end{document}